\documentclass[12pt,reqno]{amsart}
\usepackage{enumerate}
\usepackage{amsmath, amsfonts, amssymb, amsthm}
\def\pmod #1{\ ({\rm{mod}}\ #1)}
\def\Z{\Bbb Z}
\def\N{\Bbb N}

\def\l{\left}
\def\r{\right}
\def\bg{\bigg}
\def\({\bg(}
\def\){\bg)}
\def\]{]\!]}
\def\[{[\![}
\def\t{\text}
\def\f{\frac}

\def\ls{\leqslant}

\def\bi{\binom}
\def\al{\alpha}

\def\eq{\equiv}
\def\Da{\Delta}

\def\jacob #1#2{\genfrac{(}{)}{}{}{#1}{#2}}

\def\Proof{\noindent{\it Proof}}
\theoremstyle{plain}
\newtheorem{theorem}{Theorem}
\newtheorem{lemma}{Lemma}

\theoremstyle{definition}

\theoremstyle{remark}
\newtheorem{remark}{Remark}

 \vspace{4mm}

\begin{document}
 \baselineskip=17pt
\hbox{Colloq. Math. 139(2015), no.\,1, 127--136.}
\medskip

\title
[Some congruences involving binomial coefficients] {Some congruences involving binomial coefficients}

\author
[Hui-Qin Cao and Zhi-Wei Sun] {Hui-Qin Cao and Zhi-Wei Sun}

\thanks{Supported by the National Natural Science Foundation (Grant Nos. 11201233 and 11171140) of China. The second author is the corresponding author.}

\address{(Hui-Qin Cao) Department of Applied Mathematics, Nanjing Audit University,
Nanjing 211815, People's Republic of China} \email{caohq@nau.edu.cn}
\address {(Zhi-Wei Sun) Department of Mathematics, Nanjing
University, Nanjing 210093, People's Republic of China}
\email{zwsun@nju.edu.cn}

\keywords{Congruences, binomial coefficients, Lucas sequences,
central trinomial coefficients
\newline \indent 2010 {\it Mathematics Subject Classification}. Primary 11A07, 11B65; Secondary 05A10, 05A19, 11B39.}

 \begin{abstract} Binomial coefficients and central trinomial coefficients play important roles in combinatorics.
 Let $p>3$ be a prime. We show that
 $$T_{p-1}\eq\l(\f p3\r)3^{p-1}\ \pmod{p^2},$$
 where the central trinomial coefficient $T_n$ is the constant term in the expansion of $(1+x+x^{-1})^n$.
 We also prove three congruences modulo $p^3$ conjectured by Sun, one of which is 
$$\sum_{k=0}^{p-1}\bi{p-1}k\bi{2k}k((-1)^k-(-3)^{-k})\eq
\l(\f p3\r)(3^{p-1}-1)\ \pmod{p^3}.$$
In addition, we get some new combinatorial identities.
\end{abstract}

\maketitle

\section{Introduction}

Throughout this paper, we set $\N=\{0,1,2,\ldots\}$ and $\Z^+=\{1,2,3,\ldots\}$.

Let $A,B\in\Z$. The Lucas sequences $u_n=u_n(A,B)\ (n\in\N)$ and
$v_n=v_n(A,B)\ (n\in\N)$ are defined by
$$
u_0=0,\ u_1=1,\ \text{and } u_{n+1}=Au_n-Bu_{n-1}\ (n\in\Z^+)
$$
and
$$
v_0=2,\ v_1=A, \ \text{and } v_{n+1}=Av_n-Bv_{n-1}\ (n\in\Z^+).
$$
The roots of the characteristic equation $x^2-Ax+B=0$ are
$$
\alpha=\frac{A+\sqrt{\Delta}}{2}\quad\text{and}\quad\beta=\frac{A-\sqrt{\Delta}}{2},
$$
where $\Delta=A^2-4B$. By induction, one can easily deduce  the
following known formulae:
$$(\alpha-\beta)u_n=\alpha^n-\beta^n\quad\text{and}\quad
v_n=\alpha^n+\beta^n\quad\t{for any}\ n\in\N.$$
(Note that in the case $\Delta=0$ we have $v_n=2(A/2)^n$ for all $n\in\N$.) It is well-known that
\begin{equation}\label{1.1}u_p\eq\jacob{\Da}{p}\pmod{p}\quad\text{and}\quad u_{p-\jacob{\Da}{p}}\eq0\pmod{p}\end{equation}
for any odd prime $p$ not dividing $B$ (see, e.g., Sun \cite{SSci}),
where $(-)$ denotes the Legendre symbol.

Let $p>3$ be a prime and let $m$ be an integer not divisible by $p$. Recently, Sun \cite{SSci, SColl}
established the following general congruences involving central binomial coefficients
and Lucas sequences:
\begin{equation}\label{bm}
\sum_{k=0}^{p-1}\frac{\bi{2k}{k}}{m^k}\eq\jacob{\Da}{p}+u_{p-\jacob{\Da}{p}}(m-2,1)\pmod{p^2}
\end{equation}
and
\begin{equation}\label{pbm}
\sum_{k=0}^{p-1}\bi{p-1}{k}\frac{\bi{2k}{k}}{(-m)^k}\eq\jacob{\Da}{p}(m-4)^{p-1}+\l(1-\frac
m2\r)u_{p-\jacob{\Da}{p}}(m-2,1)\pmod{p^2},
\end{equation}
where $\Delta=m^2-4m$. Clearly $\bi{p-1}k\eq(-1)^k\pmod p$ for all $k=0,\ldots,p-1$.

Note that for each $n=0,1,2,\ldots$ the central binomial coefficient $\bi{2n}n$ is the constant term of $(1+x)^{2n}/x^n=(2+x+x^{-1})^n$.
For $n\in\N$, the {\it central trinomial coefficient} $T_n$ is the constant term in the expansion of $(1+x+x^{-1})^n$, i.e.,
$$T_n=\sum_{k=0}^{\lfloor n/2\rfloor}\f{n!}{k!k!(n-2k)!}=\sum_{k=0}^{\lfloor n/2\rfloor}\bi nk\bi{n-k}k.$$
Central trinomial coefficients arise naturally in enumerative combinatorics (cf. Sloane \cite{Sl}), e.g.,
$T_n$ is the number of lattice paths from the point $(0, 0)$ to $(n, 0)$ with only allowed steps $(1,0)$, $(1, 1)$ and $(1, -1)$.
As Andrews \cite{A} pointed out, central trinomial coefficients were first studied by L. Euler.
Recently, Sun \cite{SCM} investigated congruence properties of central trinomial coefficients; for example, he proved that
$\sum_{k=0}^{p-1}T_k^2\eq(\f{-1}p)\pmod p$ for any odd prime $p$.

Now we state our first theorem.

\begin{theorem}
Let $p>3$ be a prime.

{\rm (i)} We have
\begin{equation}\label{tt}
T_{p-1}\eq\l(\f p3\r)3^{p-1}\pmod{p^2}
\end{equation}
and
\begin{equation}\label{1.4}\sum_{k=0}^{p-1}\bi{p-1}{k}\bi{2k}{k}((-1)^k-(-3)^{-k})\eq\jacob{p}{3}(3^{p-1}-1)\pmod{p^3}.\end{equation}

{\rm (ii)} If $p\eq\pm 1\pmod{12}$, then
\begin{equation}\label{u41}\sum_{k=0}^{p-1}\bi{p-1}{k}\bi{2k}{k}(-1)^ku_k(4,1)\eq(-1)^{(p-1)/2}u_{p-1}(4,1)\pmod{p^3}.\end{equation}
If $p\eq\pm 1\pmod{8}$, then
\begin{equation}\label{u42}\sum_{k=0}^{p-1}\bi{p-1}{k}\bi{2k}{k}\frac{u_k(4,2)}{(-2)^k}\eq(-1)^{(p-1)/2}u_{p-1}(4,2)\pmod{p^3}.\end{equation}
\end{theorem}
\begin{remark} (\ref{1.4}) and part (ii) of Theorem 1.1 were conjectured by Sun \cite[Conj. 1.3]{SPMD}.
\end{remark}

During our efforts to prove Theorem 1.1, we also obtain some combinatorial identities.

\begin{theorem} Let $n$ be a positive integer.

{\rm (i)} If $6\ |\ n$, then
\begin{equation}\label{1u}
\sum_{k=0}^n\binom nk\binom{2k}{k}\frac{\jacob k3}{4^k}=0.
\end{equation}
If $n\equiv 3\pmod{6}$, then
\begin{equation}\label{1v}
\sum_{k=0}^n\binom nk\binom{2k}{k}\frac{3[3|k]-1}{4^k}=0,
\end{equation}
where $[3\mid k]$ is $1$ or $0$ according as $3\mid k$ or not.

{\rm (ii)} If $4\ |\ n$, then
\begin{equation}\label{2u}
\sum_{k=0}^n\binom nk\binom{2k}{k}\frac{u_k(2,2)}{(-4)^k}=0.
\end{equation}
If $n\equiv 2\pmod{4}$, then
\begin{equation}\label{2v}
\sum_{k=0}^n\binom nk\binom{2k}{k}\frac{v_k(2,2)}{(-4)^k}=0.
\end{equation}

{\rm (iii)} If $3\mid n$, then
\begin{equation}\label{3u}
\sum_{k=0}^n\binom nk\binom{2k}{k}\frac{u_k(3,3)}{(-4)^k}=0.
\end{equation}
\end{theorem}

We will provide two lemmas in the next section and prove Theorems 1.1 and 1.2 in Section 3.

\section{Two Lemmas}\setcounter{equation}{0}

\begin{lemma}\label{luv}
Let $A\in\Z^+$ and $B,m\in\Z\setminus\{0\}$ with $\Delta=A^2-4B\not=0$. Let $\al=(A+\sqrt{\Delta})/2$ and $\beta=(A-\sqrt{\Delta})/2$. Then, for every $n\in\N$ we have
\begin{equation}\label{uab}
\sum_{k=0}^n\binom
nk\binom{2k}{k}\frac{u_k(A,B)}{m^k}=\frac{d^{n/2}(\alpha^n-(-\beta)^n)}{m^n(\alpha-\beta)}\sum_{k=0}^{\lfloor
n/2\rfloor}\binom{n}{k}\bi{n-k}{k}d^{-k}
\end{equation}
and
\begin{equation}\label{vab}
\sum_{k=0}^n\binom
nk\binom{2k}{k}\frac{v_k(A,B)}{m^k}=\frac{d^{n/2}(\alpha^n+(-\beta)^n)}{m^n}\sum_{k=0}^{\lfloor
n/2\rfloor}\binom{n}{k}\bi{n-k}{k}d^{-k},
\end{equation}
where $m=-4B/A$ and $d=4\Delta/A^2$.
\end{lemma}

\Proof. For a polynomial $P(x)$ over the field of complex numbers, we use $[x^n]P(x)$ to denote the
coefficient of $x^n$ in $P(x)$. It's easy to see that
\begin{align*}
[x^n]((1+\alpha x)^2+mx)^n&=[x^n]\sum_{k=0}^n\binom nk(1+\alpha
x)^{2k}(mx)^{n-k}\\
&=m^n\sum_{k=0}^n\binom nk\binom{2k}{k}\frac{\alpha^k}{m^k}.
\end{align*}
On the other hand,
\begin{align*}
[x^n]((1+\alpha x)^2+mx)^n&=[x^n](\alpha^2x^2+(2\alpha+m)x+1)^n\\
&=[x^n]\sum_{\substack{r,s,t\geq0\\r+s+t=n}}\binom{n}{r,s,t}\alpha^{2r}(2\alpha+m)^sx^{2r+s}\\
&=\alpha^{n}\sum_{\substack{r,s\geq0\\2r+s=n}}\binom{n}{r,s,r}\(2+\frac{m}{\alpha}\)^s\\
&=\alpha^{n}\sum_{k=0}^{\lfloor
n/2\rfloor}\bi{n}{k}\bi{n-k}{k}\(2+\frac{m}{\alpha}\)^{n-2k}.
\end{align*}
So we obtain
\begin{equation}\label{alpha}
m^n\sum_{k=0}^n\binom
nk\binom{2k}{k}\frac{\alpha^k}{m^k}=\alpha^{n}\sum_{k=0}^{\lfloor
n/2\rfloor}\bi{n}{k}\bi{n-k}{k}\(2+\frac{m}{\alpha}\)^{n-2k}.
\end{equation}
Similarly,
\begin{equation}\label{beta}
m^n\sum_{k=0}^n\binom
nk\binom{2k}{k}\frac{\beta^k}{m^k}=\beta^{n}\sum_{k=0}^{\lfloor
n/2\rfloor}\bi{n}{k}\bi{n-k}{k}\(2+\frac{m}{\beta}\)^{n-2k}.
\end{equation}
As $4B=-mA$, we see that
\begin{align*}
2+\frac{2m}{A\pm\sqrt{\Delta}}=2+\frac{2m(A\mp\sqrt{\Delta})}{4B}=\pm\f{2m}{mA}\sqrt{A^2+mA}=\pm\sqrt{d},
\end{align*}
i.e., $2+m/\al=\sqrt d$ and $2+m/\beta=-\sqrt d$.
Since $u_k=(\alpha^k-\beta^k)/(\al-\beta)$ and
$v_k=\alpha^k+\beta^k$ for all $k\in\N$, combining (\ref{alpha}) and (\ref{beta}) we
get (\ref{uab}) and (\ref{vab}) immediately. \qed

\begin{lemma}
Let $p>3$ be a prime, and let $d\in\Z$ with $p\nmid d$. Then
\begin{equation}\label{2a}\begin{aligned}
&\sum_{k=0}^{(p-1)/2}\bi{p-1}{k}\bi{p-1-k}{k}d^{-k}\\
\eq&\jacob{D}{p}\(\frac{1-d^{p-1}}{2}+(d-4)^{p-1}\)-\frac{d}{4}u_{p-\jacob{D}{p}}(d-2,1)\pmod{p^2},
\end{aligned}\end{equation}
where $D=d(d-4)$.
\end{lemma}
\Proof. For every $k=0,1,\ldots,p-1$, we clearly have
\begin{equation}\label{2.6}\binom{p-1}k=(-1)^k\prod_{0<j\ls k}\l(1-\f pj\r)\eq (-1)^k(1-pH_k)\pmod{p^2},\end{equation}
where $H_k$ denotes the harmonic number $\sum_{0<j\ls k}1/j$.
Thus
\begin{align*}
&\sum_{k=0}^{(p-1)/2}\binom{p-1}{k}\binom{p-1-k}{k}d^{-k}\\
\equiv&\sum_{k=0}^{(p-1)/2}(-1)^k(1-pH_k)\binom{p-1-k}{k}d^{-k}\\
=&\sum_{k=0}^{(p-1)/2}\binom{p-1-k}{k}(-d)^{-k}-p\sum_{k=0}^{(p-1)/2}H_k\binom{p-1-k}{k}(-d)^{-k}\pmod{p^2}.
\end{align*}
Since $\bi{p-1-k}k\equiv\bi{-1-k}k=(-1)^k\binom{2k}k\pmod{p}$ for all $k=0,\ldots,p-1$, we obtain from the above
\begin{equation}\label{aa}\begin{aligned}&\sum_{k=0}^{(p-1)/2}\binom{p-1}{k}\binom{p-1-k}{k}d^{-k}
\\\equiv&\sum_{k=0}^{(p-1)/2}\binom{p-1-k}{k}(-d)^{-k}-p\sum_{k=0}^{(p-1)/2}H_k\binom{2k}{k}d^{-k}\pmod{p^2}.
\end{aligned}\end{equation}

It is known that
$$u_{n+1}(A,
B)=\sum_{k=0}^{\lfloor n/2\rfloor}\bi{n-k}{k}A^{n-2k}(-B)^k\quad\t{for all}\ \quad
n=0,1,2,\ldots$$ which can be easily proved by induction.
So we have
$$u_p(d,d)=\sum_{k=0}^{(p-1)/2}\bi{p-1-k}kd^{p-1-2k}(-d)^k=d^{p-1}\sum_{k=0}^{(p-1)/2}\bi{p-1-k}k(-d)^{-k}.$$
By \cite[Lemma 2.4]{SSci},
\begin{equation*}
2u_p(d,d)-\jacob{D}{p}d^{p-1}\eq
u_p(d-2,1)+u_{p-\jacob{D}{p}}(d-2,1)\pmod{p^2}.
\end{equation*}
In view of [4, (3.6)], if $p\nmid d-4$ then
\begin{equation*}
u_p(d-2,1)-\jacob{D}{p}\eq\(\frac
d2-1\)u_{p-\jacob{D}{p}}(d-2,1)\pmod{p^2}.
\end{equation*}
This also holds when $p\mid d-4$, since $(\f Dp)=0$ and
$$u_p(d-2,1)=u_{p-(\f Dp)}(d-2,1)=u_{p-\l(\f{(d-2)^2-4\cdot1}p\r)}(d-2,1)\eq0\ \pmod p$$
by (\ref{1.1}). Combining the above two congruences we immediately get
$$u_p(d,d)\eq\jacob{D}{p}\frac{d^{p-1}+1}{2}
+\frac{d}{4}u_{p-\jacob{D}{p}}(d-2,1)\ \pmod{p^2}.$$ Hence
\begin{equation}\label{dd}\sum_{k=0}^{(p-1)/2}\binom{p-1-k}{k}(-d)^{-k}
\eq\jacob{D}{p}\frac{d^{p-1}+1}{2d^{p-1}}+\frac{d}{4}u_{p-\jacob{D}{p}}(d-2,1)\pmod{p^2}
\end{equation}
since $u_{p-(\frac{D}p)}(d-2,1)\equiv0\pmod p$ and $d^{p-1}\equiv 1\pmod p$.

Note that $p\mid\bi{2k}k$ for $k=(p+1)/2,\ldots,p-1$. With the help of (\ref{2.6}), we have
\begin{align*}
&p\sum_{k=0}^{(p-1)/2}H_k\binom{2k}{k}d^{-k}\\
\equiv&\sum_{k=0}^{(p-1)/2}\left(1-(-1)^k\binom{p-1}{k}\right)\binom{2k}{k}d^{-k}
\\=&\sum_{k=0}^{(p-1)/2}\frac{\binom{2k}{k}}{d^k}-\sum_{k=0}^{(p-1)/2}\binom{p-1}{k}\frac{\binom{2k}{k}}{(-d)^k}\\
=&\sum_{k=0}^{(p-1)/2}\frac{\binom{2k}{k}}{d^k}+\sum_{k=(p+1)/2}^{p-1}\binom{p-1}{k}\frac{\binom{2k}{k}}{(-d)^k}-\sum_{k=0}^{p-1}\binom{p-1}{k}\frac{\binom{2k}{k}}{(-d)^k}\\
\equiv&\sum_{k=0}^{p-1}\frac{\binom{2k}{k}}{d^k}-\sum_{k=0}^{p-1}\binom{p-1}{k}\frac{\binom{2k}{k}}{(-d)^k}\ \pmod{p^2}.
\end{align*}
Thus, by applying (\ref{bm}) and (\ref{pbm}) with $m=d$ we find that $p\sum_{k=0}^{(p-1)/2}H_k\binom{2k}{k}d^{-k}$ is congruent to
$$\jacob{D}{p}+u_{p-\jacob{D}{p}}(d-2,1)-\(1-\frac{d}{2}\)u_{p-\jacob{D}{p}}(d-2,1)-\jacob{D}{p}(d-4)^{p-1}$$
modulo $p^2$. Thus

\begin{equation}\label{hh}p\sum_{k=0}^{(p-1)/2}H_k\binom{2k}{k}d^{-k}
\eq\jacob{D}{p}(1-(d-4)^{p-1})+\frac{d}{2}u_{p-\jacob{D}{p}}(d-2,1)\ \pmod{p^2}.
\end{equation}

Combining (\ref{aa}), (\ref{dd}) and (\ref{hh}), we finally obtain
\begin{align*}
&\sum_{k=0}^{(p-1)/2}\binom{p-1}{k}\binom{p-1-k}{k}d^{-k}\\
\eq&\jacob{D}{p}\(\frac{1-d^{p-1}}{2d^{p-1}}+(d-4)^{p-1}\)-\frac{d}{4}u_{p-\jacob{D}{p}}(d-2,1)\\
\eq&\jacob{D}{p}\(\frac{1-d^{p-1}}{2}+(d-4)^{p-1}\)-\frac{d}{4}u_{p-\jacob{D}{p}}(d-2,1)\pmod{p^2}.
\end{align*}
This concludes the proof. \qed

\section{Proofs of Theorems 1.1 and 1.2}
\setcounter{equation}{0}

\noindent{\it Proof of Theorem} 1.1(i).
Let $\omega$ be the primitive cubic root $(-1+\sqrt{-3})/2$. For each $k=0,1,2,\ldots$, we clearly have
$$u_{3k}(-1,1)=u_{3k}(\omega+\bar\omega,\omega\bar\omega)=\f{\omega^{3k}-\bar\omega^{3k}}{\omega-\bar\omega}=0.$$
As
$$T_{p-1}=\sum_{k=0}^{(p-1)/2}\bi{p-1}k\bi{p-1-k}k,$$
applying (\ref{2a}) with $d=1$ we get
$$T_{p-1}\eq\left(\f{-3}p\right)(-3)^{p-1}-\f14u_{p-(\f{-3}p)}(-1,1)=\l(\f p3\r)3^{p-1}\pmod{p^2}.$$
This prove (\ref{tt}).

Note that $u_k(4,3)=(3^k-1)/(3-1)$ for all $k\in\N$. With the help of Lemma 2.1 and (\ref{tt}), we have
\begin{align*}
&\sum_{k=0}^{p-1}\binom
{p-1}{k}\binom{2k}{k}\frac{u_k(4,3)}{(-3)^k}\\
=&\frac{3^{p-1}-(-1)^{p-1}}{(3-1)(-3)^{p-1}}\sum_{k=0}^{(p-1)/2}\bi{p-1}{k}\bi{p-1-k}{k}=\f{3^{p-1}-1}{2\times3^{p-1}}T_{p-1}
\\\eq&\f{3^{p-1}-1}{2\times3^{p-1}}\l(\f p3\r)3^{p-1}\pmod{p^3}
\end{align*}
and hence the desired (\ref{1.4}) follows. \qed

\noindent{\it Proof of Theorem} 1.1(ii). Suppose that $p\eq\pm 1\pmod{12}$. In light of the second congruence in (\ref{1.1}),
$$u_{p-1}(4,1)=u_{p-(\f{4^2-4\cdot1}p)}(4,1)\eq0\pmod p.$$
By Lemma 2.2,
\begin{align*}
&\sum_{k=0}^{(p-1)/2}\bi{p-1}k\bi{p-1-k}k3^{-k}
\\\eq&\l(\f{-3}p\r)\l(\f{1-3^{p-1}}2+(-1)^{p-1}\r)-\f34u_{p-(\f{-3}p)}(1,1)
\eq\l(\f p3\r)\f{3-3^{p-1}}2\pmod{p^2}
\end{align*}
since
$$u_{3k}(1,1)=\f{(-\omega)^{3k}-(-\bar\omega)^{3k}}{-\omega-(-\bar\omega)}=0\quad\mbox{for all}\ k\in\N.$$
Combining this with Lemma 2.1 we get
\begin{align*}&\sum_{k=0}^{p-1}\bi{p-1}k\bi{2k}k(-1)^ku_k(4,1)
\\=&\f{3^{(p-1)/2}}{(-1)^{p-1}}u_{p-1}(4,1)\sum_{k=0}^{(p-1)/2}\bi{p-1}k\bi{p-1-k}k3^{-k}
\\\eq&3^{(p-1)/2}u_{p-1}(4,1)\l(\f p3\r)\f{3-3^{p-1}}2\pmod{p^3}.
\end{align*}
Note that $3^{p-1}\eq
2\cdot3^{(p-1)/2}-1\pmod{p^2}$ since $3^{(p-1)/2}\eq(\f3p)=1\pmod p$. So we have
\begin{align*}
\sum_{k=0}^{p-1}\binom
{p-1}{k}\binom{2k}{k}(-1)^ku_k(4,1)
\eq&3^{(p-1)/2}\jacob{-3}{p}\frac{3-3^{p-1}}{2}u_{p-1}(4,1)
\\\eq&(-1)^{(p-1)/2}3^{(p-1)/2}(2-3^{(p-1)/2})u_{p-1}(4,1)\\
\eq&(-1)^{(p-1)/2}u_{p-1}(4,1)\ \pmod{p^3}.
\end{align*}
This proves (\ref{u41}).
\medskip

Now assume that $p\eq\pm1\pmod 8$. In view of the second congruence in (\ref{1.1}),
$$u_{p-1}(4,2)=u_{p-(\f{4^2-4\cdot2}p)}(4,2)\eq0\pmod p.$$
By Lemma 2.2,
\begin{align*}
&\sum_{k=0}^{(p-1)/2}\bi{p-1}k\bi{p-1-k}k2^{-k}
\\\eq&\l(\f{-4}p\r)\l(\f{1-2^{p-1}}2+(-2)^{p-1}\r)-\f24u_{p-(\f{-4}p)}(0,1)
=\l(\f {-1}p\r)\f{1+2^{p-1}}2\ \pmod{p^2}
\end{align*}
since $u_{2k}(0,1)=0$ for all $k\in\N$.
Combining this with Lemma 2.1 we get
\begin{align*}\sum_{k=0}^{p-1}\bi{p-1}k\bi{2k}k\f{u_k(4,2)}{(-2)^k}=&\f{2^{(p-1)/2}}{(-2)^{p-1}}u_{p-1}(4,2)\sum_{k=0}^{(p-1)/2}\bi{p-1}k\bi{p-1-k}k2^{-k}
\\\eq&\f{u_{p-1}(4,2)}{2^{(p-1)/2}}\l(\f {-1}p\r)\f{1+2^{p-1}}2\pmod{p^3}.
\end{align*}
This is equivalent to (\ref{u42}) since $2^{p-1}+1-2\cdot2^{(p-1)/2}=(2^{(p-1)/2}-1)^2\eq0\pmod{p^2}$.

In view of the above, we have completed the proof of Theorem 1.1(ii). \qed
\medskip

\noindent{\it Proof of Theorem 1.2}. (i) As $-\omega-\bar\omega=1$ and
$(-\omega)(-\bar\omega)=1$, for any $k\in\Z$ we have
$$u_k(1,1)=\f{(-\omega)^k-(-\bar\omega)^k}{-\omega-(-\bar\omega)}=(-1)^{k-1}\jacob{k}{3}$$
and
$$ v_k(1,1)=(-\omega)^k+(-\bar\omega)^k=(-1)^k(3[3|k]-1).$$
If $6\mid n$, then $(-\omega)^n=1=\bar\omega^n$ and hence by (\ref{uab}) we have
$$\sum_{k=0}^n\bi nk\bi{2k}k\f{u_k(1,1)}{(-4)^k}=0,$$
which is equivalent to (\ref{1u}). If $n\eq3\pmod 6$, then
$(-\omega)^n=-1=-\bar\omega^n$ and hence by (\ref{vab}) we have
$$\sum_{k=0}^n\bi nk\bi{2k}k\f{v_k(1,1)}{(-4)^k}=0,$$
which is equivalent to (\ref{1v}).

(ii) Clearly $(1+i)+(1-i)=(1+i)(1-i)=2$.
When $n$ is even,
$$(1+i)^n=i^n(1-i)^n=(-1)^{n/2}(1-i)^n=\begin{cases}(i-1)^n&\mbox{if}\ 4\mid n,
\\-(i-1)^n&\mbox{if}\ n\eq2\pmod4.\end{cases}$$
So we get the desired result in Theorem 1.2(ii) by applying Lemma 2.1.

(iii) Let $\alpha=(3+\sqrt{-3})/2$ and $\beta=(3-\sqrt{-3})/2$. Then $\al+\beta=\al\beta=3$.
Observe that
$$\al^2-\al\beta+\beta^2=(\al+\beta)^2-3\al\beta=0$$
and hence $\al^3=(-\beta)^3$. If $3\mid n$, then $\al^n=(-\beta)^n$ and hence
(\ref{3u}) holds by (\ref{uab}).

\medskip

In view of the above, we have finished the proof of Theorem 1.2. \qed

\medskip
\noindent
{\bf Acknowledgment}. The authors would like to thank the referee for helpful comments.


\medskip

\begin{thebibliography}{99}

\bibitem{A} G. E. Andrews, {\it Euler's ``exemplum memorabile inductionis fallacis" and $q$-trinomial coefficients},
J. Amer. Math. Soc. {\bf 3} (1990), 653--669.

\bibitem{Sl} N. J. A. Sloane, Sequence A002426 in OEIS
(On-Line Encyclopedia of Integer Sequences), {\tt http://oeis.org}.

\bibitem {SSci} Z.-W. Sun,
\textit{Binomial coefficients, Catalan numbers and Lucas quotients},
Sci. China Math. {\bf 53} (2010), 2473--2488.

\bibitem {SColl} Z.-W. Sun, \textit{On sums of binomial coefficients modulo $p^2$},
Colloq. Math. {\bf 127} (2012), 39--54.

\bibitem {SPMD} Z.-W. Sun,
\textit{On harmonic numbers and Lucas sequences},
Publ. Math. Debrecen {\bf 80} (2012), 25--41.

\bibitem {SCM} Z.-W. Sun,
\textit{Congruences involving generalized central binomial coefficients},
Sci. China Math. {\bf 57} (2014), 1375--1400.

\end{thebibliography}
\end{document}